\documentclass{amsart}

\usepackage{amssymb}

\usepackage{amsmath,amscd,graphicx,amsthm}

\newtheorem{theorem}{Theorem}

\newtheorem{proposition}[theorem]{Proposition}
\newtheorem{corollary}[theorem]{Corollary}

\theoremstyle{definition}

\theoremstyle{remark}
\newtheorem{remark}[theorem]{Remark}

\numberwithin{equation}{section}
\numberwithin{theorem}{section}

\def\A{{\mathcal A}}
\def\AA{{\mathbb A}}

\def\C{{\mathbb C}}
\def\CC{{\mathcal C}}
\def\D{{\mathfrak d}}

\def\FF{{\mathcal F}}
\def\FFF{{\mathcal F}}
\def\G{{\mathcal G}}
\def\GCC{{\G\CC}}

\def\O{{\mathcal O}}
\def\P{{\mathcal P}}

\def\TE{{\mathcal T}}

\def\UU{\overline{\A}}

\def\Z{{\mathbb Z}}

\def\fy{\varphi}
\def\g{\mathfrak g}

\def\one{\mathbf 1}

\def\r{{\bf r}}

\def\wB{{\widetilde{B}}}

\def\wx{{\widetilde{\bf x}}}
\def\x{{\bf x}}

\def\End{\operatorname{End}}

\def\Id{{\operatorname {Id}}}

\def\Mat{\operatorname{Mat}}

\def\Poi{{\{\cdot,\cdot\}}}

\def\dnabla{{\raisebox{2pt}{$\bigtriangledown$}}\negthinspace}

\begin{document}
\title{Generalized cluster structure on the Drinfeld double of $GL_n$}

\author{M. Gekhtman}

\address{Department of Mathematics, University of Notre Dame, Notre Dame,
IN 46556}
\email{mgekhtma@nd.edu}

\author{M. Shapiro}
\address{Department of Mathematics, Michigan State University, East Lansing,
MI 48823}
\email{mshapiro@math.msu.edu}

\author{A. Vainshtein}
\address{Department of Mathematics \& Department of Computer Science, University of Haifa, Haifa,
Mount Carmel 31905, Israel}
\email{alek@cs.haifa.ac.il}

\begin{abstract}
We construct a generalized cluster structure compatible with  the Poisson bracket on the Drinfeld double of the standard
Poisson-Lie group $GL_n$ and derive from it a generalized cluster structure in $GL_n$ compatible with the push-forward of the dual Poisson--Lie bracket.

\end{abstract}
\maketitle

\section{Introduction}
\label{intro}

The connection between cluster algebras and Poisson structures is documented in \cite{GSVb}. Among the most important examples in which this connection has been utilized are coordinate rings of double Bruhat cells in semisimple Lie groups equipped with (the restriction of) the standard Poisson-Lie structure. In \cite{GSVb}, we applied our technique of constructing a cluster structure compatible with a given Poisson structure in this situation and recovered the cluster structure built in \cite{CAIII}. The standard Poisson-Lie structure is a particular case of Poisson-Lie structures corresponding to quasi-triangular Lie bialgebras. Such structures are associated with solutions to the classical Yang-Baxter equation. Their complete classification was obtained by Belavin and Drinfeld in \cite{BD} in terms of certain combinatorial data defined in terms of the corresponding root system. In \cite{GSVMMJ} we conjectured that any such solution gives rise to a compatible cluster structure on the Lie group and provided several examples supporting this conjecture. 
%by showing that it holds true for the class of the standard Poisson-Lie structure in any simple complex Lie group, and for the whole Belavin-Drinfeld classification in $GL_n$ for $n = 2, 3, 4$. (The $GL_5$ case was treated in \cite{Idan}.) 
%We call the cluster structures associated with the nontrivial Belavin-Drinfeld data {\em exotic}. 
Recently \cite{GSVPNAS,GSVMem}, we constructed the 
%exotic 
cluster structure corresponding to the Cremmer--Gervais Poisson structure in $GL_n$ for any $n$.

As we established in \cite{GSVMem}, the construction of 
%exotic 
cluster structures on a simple Poisson-Lie group $\G$  relies on properties of the Drinfeld double $D(\G)$. Moreover, in the Cremmer-Gervais case generalized determinantal identities on which cluster transformations are modeled can be extended to identities valid in the double. It is not too far-fetched then to suspect that 
%in the double 
there exists a cluster structure on $D(\G)$ compatible with the Poisson-Lie bracket induced by the Poisson-Lie bracket on $\G$.
However, an interesting phenomenon was observed even in the first nontrivial example of $D(GL_2)$: although we were able to construct a log-canonical regular coordinate chart in terms of which all standard coordinate functions are expressed as (subtraction free) Laurent polynomials, it is not possible to define cluster transformations in such a way that all cluster variables that one expects to be mutable transform into regular functions. This problem is resolved, however, if one is allowed to use {\em generalized cluster transformations} previously considered in \cite{GSV1,GSVb} and, more recently, axiomatized in  \cite{CheSha}.

%This%version of cluster transformations, though not as general is the one used in the definition%of LP-algebras, has an advantage of still being amenable to treatment by tools of Poisson%geometry. For example, D(SL2) does admit a generalized cluster structure of finite type%in the sense of [9] ( the corresponding exchange graph is isomorphic to the 3-dimensional%Bott-Taubes polytop).{\em We conjecture that for any quasi-triangular Poisson-Lie group $(\G, \{\ , \ \})$, (i) the ring ofregular functions on $D(\G)$ admits a generalized cluster structure compatible with the correspondingPoisson-Lie structure on $D(\G)$ and (ii) this generalized cluster structure restricts tothe exotic cluster structure on $\G$ compatible with $\Poi$ and to a generalized cluster structure onthe dual Poisson-Li group $\G$.}

In this note, we describe
such a generalized cluster structure on the Drinfeld double in the case of the standard Poisson-Lie group $GL_n$.
Using this structure, one can recover the standard cluster structure on $GL_n$  and introduce a generalized cluster structure on $GL_n$ compatible with the Poisson
 bracket dual to the standard Poisson--Lie bracket. Note that the log-canonical basis suggested in \cite{Bra} is different from the one
constructed here and does not lead to a regular cluster structure.

\section{Generalized cluster structures of geometric type and compatible Poisson brackets}
\label{SecPrel}

%We start with the basics on cluster algebras of geometric type. The definition that we present
%below is not the most general one, see, e.g.,
%\cite{FZ2, CAIII} for a detailed exposition. In what follows, we will use a notation $[i,j]$ for an interval
%$\{i, i+1, \ldots , j\}$ in $\mathbb{N}$, and we will denote $[1, n]$ by $[n]$.
 
%The {\em coefficient group\/} $\PP$ is a free multiplicative abelian
%group of finite rank $m$ with generators $g_1,\dots, g_m$.
%An {\em ambient field\/}  is
%the field $\FFF$ of rational functions in $n$ independent variables with
%coefficients in the field of fractions of the integer group ring
%$\Z\PP=\Z[g_1^{\pm1},\dots,g_m^{\pm1}]$ (here we write
%$x^{\pm1}$ instead of $x,x^{-1}$).

Let $\widetilde{B}=(b_{ij})$ be an $n\times(n+m)$ integer matrix
whose principal part $B$ is skew-symmetrizable (recall that the principal part of a rectangular matrix  
is its maximal leading square submatrix). 
Let $\FFF$ be the field of rational functions in $n+m$ independent variables
with rational coefficients. There are $m$  distinguished variables; we denote them
$x_{n+1},\dots,x_{n+m}$ and call {\em stable}. Finally, we define $2n$ {\em stable $\tau$-monomials\/}
$v_{i;>}$ and $v_{i;<}$, $1\le i\le n$, via 
%v_{i;>}=\prod_{\substack{n+1\le j\le n+m\\  b_{ij}>0}}x_i^{b_{ij}},\qquad
%v_{i;<}=\prod_{\substack{n+1\le j\le n+m\\  b_{ij}<0}}x_i^{-b_{ij}};
$v_{i;>}=\prod\{x_j^{b_{ij}}: n+1\le j\le n+m,  b_{ij}>0\}$, $v_{i;<}=\prod\{x_j^{-b_{ij}}: n+1\le j\le n+m,  b_{ij}<0\};$ 
here, as usual, the product over the empty set is assumed to be
equal to~$1$.

A {\em seed\/} (of {\em geometric type\/}) in $\FFF$ is a triple
$\Sigma=(\x,\widetilde{B},\P)$,
where $\x=(x_1,\dots,x_n)$ is a transcendence basis of $\FFF$ over the field of
fractions of  $\bar\AA=\Z[x_{n+1}^{\pm1},\dots,x_{n+m}^{\pm1}]$ (here we write
$x^{\pm1}$ instead of $x,x^{-1}$), and $\P$ is a set of $n$ {\em strings}. The $i$th string is a collection of 
%(Laurent monomials in stable variables?) 
polynomials $p_{ir}\in\bar\AA$, $0\le r\le d_i$, such that  $d_i$ is a factor 
of $\gcd\{b_{ij}: 1\le j\le n\}$, $p_{i0}=p_{id_i}=1$, and $\hat p_{ir}=\left(p_{ir}v_{i;>}^rv_{i;<}^{d_i-r}\right)^{1/d_i}$
belong to the polynomial ring $\AA=\Z[x_{n+1},\dots,x_{n+m}]$, $0\le r\le d_i$; it is called {\em trivial\/} if $d_i=1$, and hence both elements of the string are equal to one.

 Matrices $B$ and $\wB$ are called the
{\it exchange matrix\/} and the {\it extended exchange matrix}, respectively. The $n$-tuple  $\x$ is called a {\em cluster\/}, and its elements
$x_1,\dots,x_n$ are called {\em cluster variables\/}. The polynomials $p_{ir}$ are called {\em coefficients}.
We say that
$\widetilde{\x}=(x_1,\dots,x_{n+m})$ is an {\em extended
cluster\/}, and $\widetilde\Sigma=(\wx,\widetilde{B},\P)$ is an {\em extended seed}.

In certain cases it is convenient to represent the data $(\wB, d_1,\dots,d_n)$ by a quiver. 
Assume that 
%$(d_1^{-1},\dots,d_n^{-1})$ is a skew-symmetrizer for $B$. Then 
the matrix obtained from 
$\wB$ by replacing
each $b_{ij}$ by $b_{ij}/d_i$ for $1\le j\le n$ and retaining it for $n+1\le j\le n+m$ has a skew-symmetric 
principal part. We say that the corresponding quiver $Q$ represents $(\wB, d_1,\dots,d_n)$ and write
$\Sigma=(\x,Q,\P)$. A vertex with $d_i\ne 1$ is called
{\em special\/}, and $d_i$ is said to be its {\em order}. A stable 
%(frozen?) 
vertex $j$ such that $b_{ij}=0$, $1\le i\le n$, is called {\em isolated}. 

Given a seed as above, the {\em adjacent cluster\/} in direction $k$, $1\le k\le n$,
is defined by
$\x'=(\x\setminus\{x_k\})\cup\{x'_k\}$,
where the new cluster variable $x'_k$ is given by the {\em generalized exchange relation}
\begin{equation*}\label{exchange}
x_kx'_k=\sum_{j=0}^{d_k}\hat p_{kj}u_{k;>}^j u_{k;<}^{d_k-j}
\end{equation*}
with {\em cluster $\tau$-monomials\/} $u_{k;>}$ and $u_{k;<}$ defined by
%u_{k;>}=\prod_{\substack{1\le i\le n\\  b_{ki}>0}}x_i^{b_{ki}/d_k},\qquad
%u_{k;<}=\prod_{\substack{1\le i\le n\\  b_{ki}<0}}x_i^{-b_{ki}/d_k}.
$u_{k;>}=\prod\{x_i^{b_{ki}/d_k}: 1\le i\le n, b_{ki}>0\}$, $u_{k;<}=\prod\{x_i^{-b_{ki}/d_k}: 1\le i\le n, b_{ki}<0\}$.

We say that $\wB'$ is
obtained from $\wB$ by a {\em matrix mutation\/} in direction $k$
%and write $\wB'=\mu_k(\wB)$ 
 if
\begin{equation*}
\label{eq:MatrixMutation}
b'_{ij}=\begin{cases}
         -b_{ij}, & \text{if $i=k$ or $j=k$;}\\
                 b_{ij}+\displaystyle\frac{|b_{ik}|b_{kj}+b_{ik}|b_{kj}|}2,
                                                  &\text{otherwise.}
        \end{cases}
\end{equation*}
Note that $\gcd\{b_{ij}:1\le j\le n\}=\gcd\{b'_{ij}:1\le j\le n\}$, and for its arbitrary factor $d$, $b_{ij}=b'_{ij}\bmod d$ 
for $n+1\le j\le n+m$.

%To define {\em coefficient mutation\/} we introduce $\hat p_{ij}=\frac{p_{ij}}{v_{i;>}^{j/d_i}v_{i;<}^{1-j/d_i}}$; note that
%$\hat p_{ij}$ are, in general, Puiseux monomials (polynomials?) in stable variables, but $\hat p_{ij}^{d_i}$ are Laurent monomials (polynomials).  
The {\em coefficient mutation\/} in direction $k$ is given by
\begin{equation*}
\label{eq:CoefMutation}
 p'_{ir}=\begin{cases}
          p_{i,d_i-r}, & \text{if $i=k$;}\\
           p_{ir}, &\text{otherwise.}
        \end{cases}
\end{equation*}

Given a seed $\Sigma=(\x,\widetilde{B},\P)$, we say that a seed
$\Sigma'=(\x',\widetilde{B}',\P')$ is {\em adjacent\/} to $\Sigma$ (in direction
$k$) if $\x'$, $\widetilde{B}'$ and $\P'$ are as above. 
%is adjacent to $\x$ in direction $k$ and $\widetilde{B}'=\mu_k(\widetilde{B})$. 
Two seeds are {\em mutation equivalent\/} if they can
be connected by a sequence of pairwise adjacent seeds. 
The set of all seeds mutation equivalent to $\Sigma$ is called the {\it generalized cluster structure\/} 
(of geometric type) in $\FFF$ associated with $\Sigma$ and denoted by $\GCC(\Sigma)$; in what follows, 
we usually write $\GCC(\wB,\P)$, or even just $\GCC$ instead. Clearly, by taking $d_i=1$ for $1\le i\le n$, 
and hence making all strings trivial, we get an ordinary cluster structure.

%Following \cite{FZ2,CAIII}, 
We associate
with $\GCC(\wB,\P)$ two algebras of rank $n$ over the ground ring $\AA$:
the {\em generalized cluster algebra\/} $\A=\A(\GCC)=\A(\wB,\P)$, which 
is the $\AA$-subalgebra of $\FF$ generated by all cluster
variables in all seeds in $\GCC(\wB,\P)$, and the {\it generalized upper cluster algebra\/}
$\UU=\UU(\GCC)=\UU(\wB,\P)$, which is the intersection of the rings of Laurent polynomials over $\AA$ in cluster variables
taken over all seeds in $\GCC(\wB,\P)$. The generalized  {\it Laurent phenomenon\/} \cite{CheSha}
claims the inclusion $\A(\GCC)\subseteq\UU(\GCC)$. 
%The natural choice of the ground ring for the geometric type
%is the polynomial ring in stable variables $\AA=\Z\P_+=\Z[x_{n+1},\dots,x_{n+m}]$; this choice is assumed unless
%explicitly stated otherwise. 

Let $V$ be a quasi-affine variety over $\C$, $\C(V)$ be the field of rational functions on $V$, and
$\O(V)$ be the ring of regular functions on $V$. Let $\GCC$ be a generalized cluster structure in $\FF$ as above.
Assume that $\{f_1,\dots,f_{n+m}\}$ is a transcendence basis of $\C(V)$. Then the map $\theta: x_i\mapsto f_i$,
$1\le i\le n+m$, can be extended to a field isomorphism $\theta: \FF_\C\to \C(V)$,  
where $\FF_\C=\FF\otimes\C$ is obtained from $\FF$ by extension of scalars.
The pair $(\GCC,\theta)$ is called a generalized cluster structure {\it in\/}
$\C(V)$, $\{f_1,\dots,f_{n+m}\}$ is called an extended cluster in
 $(\GCC,\theta)$.
Sometimes we omit direct indication of $\theta$ and say that $\GCC$ is a generalized cluster structure {\em on\/} $V$. 
A generalized cluster structure $(\GCC,\theta)$ is called {\it regular\/}
if $\theta(x)$ is a regular function for any cluster variable $x$. 
The two algebras defined above have their counterparts in $\FF_\C$ obtained by extension of scalars; they are
denoted $\A_\C$ and $\UU_\C$.
If, moreover, the field isomorphism $\theta$ can be restricted to an isomorphism of 
$\A_\C$ (or $\UU_\C$) and $\O(V)$, we say that 
$\A_\C$ (or $\UU_\C$) is {\it naturally isomorphic\/} to $\O(V)$.

Let $\Poi$ be a Poisson bracket on the ambient field $\FFF$, and $\GCC$ be a generalized cluster structure in $\FFF$. 
We say that the bracket and the generalized cluster structure are {\em compatible\/} if any extended
cluster $\widetilde{\x}=(x_1,\dots,x_{n+m})$ is {\em log-canonical\/} with respect to $\Poi$, that is,
%\begin{equation}\label{cpt}
$\{x_i,x_j\}=\omega_{ij} x_ix_j$,
%\end{equation}
where $\omega_{ij}\in\Z$ are
constants for all $i,j$, $1\le i,j\le n+m$; 
%and, additionally, 
it follows that all polynomials $p_{ir}$ are Casimirs of the bracket. 
%The matrix $\Omega^{\widetilde \x}=(\omega_{ij})$ is called the {\it coefficient matrix\/}
%of $\Poi$ (in the basis $\widetilde \x$); clearly, $\Omega^{\widetilde \x}$ is skew-symmetric. 
The notion of compatibility  extends to Poisson brackets on $\FF_\C$ without any changes.

\section{Standard Poisson-Lie group $\G$ and its Drinfeld double}
\label{double}

Let $\G$ be a reductive complex Lie group equipped with a Poisson bracket $\Poi$.
$\G$ is called a {\em Poisson--Lie group\/}
if the multiplication map
$
\G\times \G \ni (x,y) \mapsto x y \in \G
$
is Poisson. 
%The most important class of Poisson--Lie groups
%is the one associated with classical R-matrices that were classified, up to an automorphism,  by Belavin and Drinfeld in \cite{BD}.
%Although this classification is important to our conjecture formulated in \cite{GSVMMJ}, we will not review it here since
%In this note
%we only deal with the case of $\G=GL_n$ equipped with the {\em standard Poisson-Lie structure} that can be described as follows.
%
%The tangent Lie algebra $\g$ of a Poisson-Lie group $\G$ has a natural structure of a {\em Lie bialgebra}. We are interested in the case when
%$\G$ be a simple complex Lie group and its tangent Lie bialgebra is {\em factorizable}. 
Denote by 
$\langle \ , \ \rangle$ an invariant nondegenerate form on
%Killing form (trace-form) on 
%$\g=sl_n$, 
$\g$, and by $\nabla^R$, $\nabla^L$ the right and
left gradients of functions on $\G$ with respect to this form. Let $\pi_{>0}, \pi_{<0}$ be projections of  
$\g$ onto subalgebras spanned by positive and negative roots and let $R=\pi_{>0} - \pi_{<0}$.
%, $\pi_0$ be the projection onto the Cartan 
%subalgebra $\h$. Define also $\pi_{\geq 0}=\pi_{>0} + \pi_0$ and $\pi_{\leq 0}=\pi_{<0} + \pi_0$.
The {\em standard Poisson-Lie bracket\/} $\Poi_r$ on $\G$  can be written as
\begin{equation}
\{f_1,f_2\}_r = \frac12\left(\langle R(\nabla^L f_1), \nabla^L f_2 \rangle - \langle R(\nabla^R f_1), \nabla^R f_2 \rangle\right).
\label{sklyabra}
\end{equation}

%Let $R_\pm\in \End\g$ be defined by $R_\pm=\frac{1}{2} ( R \pm \mbox{Id})$. 
%The images of $\g$ under $R_\pm$ are Lie subalgebras of $\g$; we denote them by $\g_\pm$,
%and the corresponding Lie subgroups of $\G$, by $\G_\pm$.
Following \cite{r-sts}, let us recall the construction of {\em the Drinfeld double}. The double of $\g$ is 
$D(\g)=\g  \oplus \g$ equipped with an invariant nondegenerate bilinear form
$\langle\langle (\xi,\eta), (\xi',\eta')\rangle\rangle = \langle \xi, \xi'\rangle - \langle \eta, \eta'\rangle$. 
Define subalgebras $\D_\pm$ of $D(\g)$ by
%\begin{equation}\label{ddeco}
$\D_+=\{( \xi,\xi) : \xi \in\g\}$ and $\D_-=\{ (R_+(\xi),R_-(\xi)) : \xi \in\g\}$,
where $R_\pm\in \End\g$ is given by $R_\pm=\frac{1}{2} ( R \pm \Id)$. 
%\end{equation}
%Then $\D_\pm$ are isotropic subalgebras of $D(\g)$ and $D(\g)= \D_+ \dot + \D_-$. In other words,
%$(D(\g), \D_+, \D_-)$ is {\em a Manin triple}. Then 
The operator $R_D= \pi_{\D_+} - \pi_{\D_-}$ can be used to define 
a Poisson--Lie structure on $D(\G)=\G\times \G$, the double of the group $\G$, via
\begin{equation}
\{f_1,f_2\}_D = \frac{1}{2}\left (\langle\langle R_D(\dnabla^L f_1), \dnabla{^L} f_2 \rangle\rangle 
- \langle\langle R_D(\dnabla^R f_1), \dnabla^R f_2 \rangle\rangle \right),
\label{sklyadouble}
\end{equation}
where $\dnabla^R$ and $\dnabla^L$ are right and left gradients with respect to $\langle\langle \cdot ,\cdot \rangle\rangle$.
Restriction of this bracket to $\G$ identified with the diagonal subgroup $\{ (X,X)\ : \ X\in \G\}$ of $D(\G)$ (whose Lie algebra is $\D_+$) 
coincides with the Poisson--Lie bracket \eqref{sklyabra}
%$\Poi_r$ 
on $\G$.

A group $\G_r$  whose Lie algebra is $\D_-$ is a Poisson-Lie subgroup of $D(\G)$ called {\em the dual Poisson-Lie group of $\G$}.
The map $\G_r \ni (X,Y) \mapsto U=X^{-1} Y$ induces another Poisson bracket on $\G$. We denote this bracket $\Poi_*$ and refer to the Poisson manifold 
$(\G, \Poi_*)$ as $\G^*$.

\section{Log-canonical basis}
\label{logcan}

In this note we only deal with the case of $\G=GL_n$.
Let $(X,Y)$ be a point in the double $D(GL_n)$. For $k,l\ge 1$, $k+l\le n-1$ define a $(k+l)\times(k+l)$ matrix 
$$
F_{kl}=F_{kl}(X,Y)=\left[\begin{array}{cc}X^{[n-k+1,n]} & Y^{[n-l+1,n]}\end{array}\right]_{[n-k-l+1,n]}.
$$ 
 For $1\le j\le i\le n$ define an $(n-i+1)\times (n-i+1)$ matrix
$$
G_{ij}=G_{ij}(X)=X_{[i,n]}^{[j,j+n-i]}.
$$
For $1\le i\le j\le n$ define an $(n-j+1)\times (n-j+1)$ matrix
$$
H_{ij}=H_{ij}(Y)=Y_{[i,i+n-j]}^{[j,n]}.
$$
For  $k,l\ge 1$, $k+l\le n$ define an $n\times n$ matrix
\begin{equation*}\label{phidef}
\Phi_{kl}=\left[\begin{array}{ccccc}(U^0)^{[n-k+1,n]}& U^{[n-l+1,n]} & (U^2)^{[n]} & \dots & (U^{n-k-l+1})^{[n]}\end{array}\right] 
\end{equation*}
where $U=X^{-1}Y$. Note that the definition of $F_{kl}$ can be extended to the case $k+l=n$ yielding $F_{n-l,l}=X\Phi_{n-l,l}$.

Denote $f_{kl}=\det F_{kl}$, $g_{ij}=\det G_{ij}$, $h_{ij}=\det H_{ij}$ and  
$$
\fy_{kl}=s_{kl}(\det X)^{n-k-l+1}\det\Phi_{kl}, 
$$
$2n^2-n+1$ functions in total. 
Here $s_{kl}$ is a sign defined as follows: it is periodic in $k+l$ with period~4 for 
$n$ odd and period~2 for $n$ even; $s_{n-l,l}=1$; $s_{n-l-1,l}=(-1)^l$ for $n$ odd and $s_{n-l-1,l}=(-1)^{l+1}$ for $n$ even;
$s_{n-l-2,l}=-1$ for $n$ odd; $s_{n-l-3,l}=(-1)^{l+1}$ for $n$ odd.
Note that the pre-factor in the definition of $\fy_{kl}$ is needed to obtain an irreducible regular function in matrix entries of $X$ and $Y$.

% We also will need a sign ${s_{kl}}$; this sign is periodic 
%in $k+l$ with period~4 for 
%$n$ odd and period~2 for $n$ even; $s_{n-l,l}=1$; $s_{n-l-1,l}=(-1)^l$ for $n$ odd and $s_{n-l-1,l}=(-1)^{l+1}$ for $n$ even;
%$s_{n-l-2,l}=-1$ for $n$ odd; $s_{n-l-3,l}=(-1)^{l+1}$ for $n$ odd; {\bf apparently, these signs can be achieved by taking the columns 
%of $\Phi_{kl}$ in the opposite order}.

Consider the polynomial $\det(X+\lambda Y)=\sum_{i=0}^n \lambda^i s_i c_i(X,Y)$, where $s_i=(-1)^i$ if $n$ is even and $s_i=1$ if $n$ is odd. It is well-known that functions $c_i(X,Y)$ are Casimirs for the Poisson-Lie bracket \eqref{sklyadouble} on $D(GL_n)$.
Note also that $c_0(X,Y)=\det X=g_{11}$ and  $c_n(X,Y)=\det Y=h_{11}$.

\begin{theorem}
\label{basis}
The family of functions $F_n=\{g_{ij}, h_{ij}, f_{kl}, \fy_{kl}, c_1,\ldots, c_{n-1} \}$ forms a log-canonical coordinate system with respect to the Poisson-Lie bracket  \eqref{sklyadouble} on $D(GL_n)$.
\end{theorem}

The proof exploits various invariance properties of functions in $F_n$.

\section{Initial quiver}
\label{init}

The quiver $Q_n$ contains $2n^2 - n +1$  vertices  labeled by the functions $g_{ij}, h_{ij}, f_{kl}, \fy_{kl}$ in the log-canonical basis $F_n$. The Casimir functions $c_1, \ldots, c_{n-1}$ correspond to isolated vertices. 
%and are not attached to vertices of $Q_n$. 
The vertex $\fy_{11}$ is special, and its order equals $n$.
The vertices $g_{i1}$, $1\le i\le n$, and $h_{1j}$, $1\le j\le n$, are stable.
%``invisible'': they are not connected to any of the other vertices, and will be not shown on figures.

The edges of $Q_n$ are comprised of $(n-1)(n-2)/2$ edge-disjoint triangles $h_{ij}\to h_{i+1,j+1}\to h_{i+1,j}\to h_{ij}$, $1\le i<j\le n-1$;
$n(n-1)/2$ disjoint triangles $g_{ij}\to g_{i+1,j+1}\to g_{i,j+1}\to g_{ij}$, $1\le j\le i\le n-1$; $(n-2)(n-3)/2$ edge-disjoint triangles
$f_{kl}\to f_{k-1,l}\to f_{k-1,l+1}\to f_{kl}$, $k+l\le n-1$, $k\ge 2$, $l\ge 1$; $(n-2)(n-3)/2$ edge-disjoint triangles
$\fy_{kl}\to \fy_{k-1,l+1}\to \fy_{k,l+1}\to \fy_{kl}$, $k+l\le n-1$, $k\ge 2$, $l\ge 1$; the path $g_{11}\to\fy_{11}\to\fy_{21}\to \fy_{12}\to \fy_{31}\to
\dots\to\fy_{1,l-1}\to\fy_{l1}\to\fy_{1l}\to\dots\to\fy_{1,n-1}$ of length $2n-3$ for $n>2$; the path $\fy_{1,n-1}\to\fy_{1,n-2}\to\dots\to\fy_{11}\to h_{11}$
of length $n-1$ for $n>2$; the path $\fy_{n-1,1}\to f_{n-2,1}\to\fy_{n-2,2}\to\dots\to\fy_{kl}\to f_{k-1,l}\to\fy_{k-1,l+1}\to\dots\to\fy_{1,n-1}$ of length
$2(n-2)$; the path $h_{11}\to\fy_{1,n-1}\to h_{22}\to f_{1,n-2}\to\dots\to f_{1l}\to h_{n-l+1,n-l+1}\to f_{1,l-1}\to\dots\to h_{nn}$ of length
$2(n-1)$; the path $h_{nn}\to h_{n-1,n}\to\dots\to h_{1n}$ of length $n-1$; the path $h_{nn}\to g_{nn}\to g_{n,n-1}\to g_{n,n-2}\to\dots\to g_{n1}$
of length $n$. Above we identify $g_{i,i+1}$ with $f_{n-i,1}$ for $1\le i\le n-1$ and $f_{k,n-k}$ with $\fy_{k,n-k}$ for $1\le k\le n-1$. 
Note that the triangle $\fy_{21}\to \fy_{12}\to \fy_{22}\to \fy_{21}$ and the first path above have a common edge; this means that for any $n>3$ there are two
edges pointing from $\fy_{21}$ to $\fy_{12}$.

\begin{figure}[ht]
\begin{center}
\includegraphics[width=12cm]{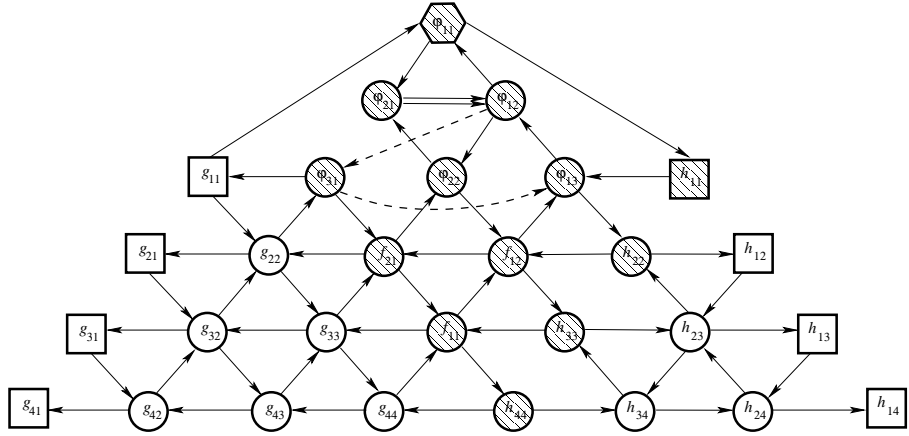}
\caption{Quiver $Q_4$}
\label{D4}
\end{center}
\end{figure}

The quiver $Q_4$ is shown in Fig.~\ref{D4}. The stable vertices are shown as squares, the special vertex is shown as a hexagon, isolated vertices are not shown. It is easy to see that $Q_4$, as well as $Q_n$ for any $n$, can be
embedded into a torus. 

\begin{remark}
{\rm
\label{diagonal} On the diagonal subgroup $\{ (X,X) : X\in GL_n\}$ of $D(GL_n)$, $g_{ii} = h_{ii}$ for $1\le i\le n$, and functions $f_{kl}$ and $\fy_{kl}$ vanish identically.
Accordingly, vertices in $Q_n$ that correspond to $f_{kl}$ and $\fy_{kl}$ are erased and, for $1\le i\le n$, vertices corresponding to $g_{ii}$ and  $h_{ii}$ are identified.
As a result, one recovers a seed of the cluster structure compatible with the standard Poisson-Lie structure on $GL_n$,  see \cite[Chap.~4.3]{GSVb}.
}

\end{remark}
 
%\begin{proposition} $Q_n$ is a planar graph.
%\end{proposition}

\section{Generalized exchange relation}
\label{long}

\begin{proposition}
\label{long_identity}
Let $A$ be a complex $n\times n$ matrix. For  $u,v\in \C^n$, define matrices
\begin{gather*}\nonumber
\Gamma(u)=\left[ u \; A u \; A^2 u \dots  A^{n-1} u\right],\\ 
\Gamma_1(u,v)=\left[ v \; u  \; A u \dots A^{n-2} u\right],\;\; 
\Gamma_2(u,v)=\left[A v \; u  \;A u  \dots  A^{n-2} u\right].
\end{gather*}
In addition, let $w$ be the last row of the classical adjoint of $\Gamma_1(u,v)$, i.e.
$w \Gamma_1(u,v) = \left (\det \Gamma_1(u,v) \right )e_n^T$. Define
$\Gamma^*(u,v)$ to be the matrix with rows $w, w A, \ldots, w A^{n-1}$. Then
\begin{equation*}
\label{longid}
\det\Big(\det\Gamma_1(u,v) A -  \det\Gamma_2(u,v) \one \Big ) = (-1)^{\frac{n(n-1)}{2}} \det\Gamma(u) \det\Gamma^*(u,v).
\end{equation*}
\end{proposition}

Specializing Proposition \ref{long_identity} to the case $A=X^{-1} Y, u=e_n, v= e_{n-1}$, one obtains

\begin{corollary}\label{polyrel}
For any $n>2$,
\begin{equation}\label{general}
\det (s_{12}\fy_{12}X+s_{21}\fy_{21}Y)=\fy_{11}P^*_n,
\end{equation}
where $P^*_n$ is a polynomial in the entries of $X$ and $Y$.
\end{corollary}

Relation \eqref{general} will serve as a generalized exchange relation in our definition of a generalized cluster structure on $D(GL_n)$. More exactly, the set
 $\P_n$ contains only one nontrivial string $\{p_{ir}\}$,
$1\le r\le n-1$. It corresponds to the vertex $\fy_{11}$, and $p_{ir}=c_r^ng_{11}^{r-n}h_{11}^{-r}$, $1\le r\le n-1$. The strings corresponding to all other vertices are trivial.

\section{Main results}
\label{main}

\begin{theorem}
\label{structure}
{\rm (i)} %The family of functions $F_n$, the matrix $\wB_n$ corresponding to the quiver $Q_n$, and the collection of coefficients $\P_n$
%define together 
%The initial generalized seed consisting of the initial quiver $Q_n$ and the family of functions $F$ that comprises both the elements of the initial extended cluster and coefficients of the together the generalized exchange relation~\eqref{general} define 
%Let $\wB_n$ be the matrix represented by the quiver $Q_n$. 
The extended seed $\widetilde\Sigma_n=(F_n,Q_n,\P_n)$ defines 
a generalized cluster structure in the ring of regular functions on $D(GL_n)$ compatible with the standard Poisson--Lie structure on $D(GL_n)$. 

{\rm (ii)} The corresponding generalized upper cluster algebra 
%coincides with 
is naturally isomorphic to
the ring of regular functions on $\Mat_n\times \Mat_n$.
\end{theorem}

Using  Theorem \ref{structure}, we can construct a  generalized cluster
structure on $GL_n^*$. For $U\in GL_n^*$, denote
$\psi_{kl}(U) = s_{kl}\det \Phi_{kl}$, where $s_{kl}$ are the signs defined in Section \ref{logcan}.
The initial extended cluster $F^*_n$  for $GL_n^*$ consists of
functions  $\psi_{kl}(U)$, $k,l\ge 1$, $k+l\le n$, $h_{ij}(U)$, $2\le i\le j\le n$, and $c_i(\one, U)$, $1\le i\le n-1$.
%, with functions $h_{ii}(U)$ serving as stable variables.
%These functions are attached to the nodes of the initial quiver $Q^*_n$.
To obtain the initial seed for $GL_n^*$, we apply a certain
sequence $\mathcal T$ of cluster transformations to the initial seed
for $D(GL_n)$. This sequence does not involve vertices associated with
functions $\fy_{kl}$. The resulting cluster $\TE(F_n)$ contains a
subset $\{ \left (\det X\right )^{\nu(f)}f  : f \in F^*_n\}$ with  $\nu(f)\in\Z_+$ (in particular,
$\nu(\psi_{kl})=n-k-l+1$).   These functions 
are attached to a subquiver $Q^*_n$ in the resulting quiver $\TE(Q_n)$, which 
is isomorphic to the subquiver of $Q_n$ formed by vertices
associated with functions $\fy_{kl}, f_{ij}$ and $h_{ii}$, see Fig.~\ref{D4}, where the vertices of the corresponding subquiver are shaded.
%The only vertices in $Q^*_n$ connected with the rest of vertices in $\TE(Q_n)$ are those associated with 
Functions $h_{ii}(U)$ are declared stable variables, $c_i(1,U)$ remain isolated.  All exchange relations  defined by mutable
vertices of $Q^*_n$ are homogeneous in $\det X$. This allows us to use $(F^*_n, Q_n^*,\P_n)$ as an initial seed for $GL_n^*$.
The generalized exchanged relation associated with the cluster variable
$\psi_{11}$ now takes form
%\begin{equation}\label{dualgeneral}
%\det (\fy_{12}\one_n+\fy_{21}U)=\fy_{11}P_n,
$\det (s_{12}\psi_{12}\one+s_{21}\psi_{21}U)=\psi_{11}\Pi^*_n$,
%\end{equation}
where $\Pi^*_n$ is a polynomial in the entries of $U$.

\begin{theorem}
\label{dualstructure}
{\rm (i)} The generalized cluster structure on $GL_n^*$ with the initial seed described above is compatible with $\Poi_*$ and regular.
 
{\rm (ii)} The corresponding generalized  upper cluster algebra 
%coincides with 
is naturally isomorphic to
the ring of regular functions on $\Mat_n$.
\end{theorem}

\section*{Acknowledgments}

M.~G.~was supported in part by NSF Grant DMS \#1362801. 
M.~S.~was supported in part by NSF Grants DMS \#1362352.  
A.~V.~was supported in part by ISF Grant \#162/12. He is grateful to Max-Planck-Institut
f\"ur Matematik, Bonn, for hospitality in September--October 2014.

\end{document}